\documentclass[12pt]{amsart}

\newtheorem{theorem}{Theorem}[section]
\newtheorem{proposition}[theorem]{Proposition}
\newtheorem{lemma}[theorem]{Lemma}
\newtheorem{remark}[theorem]{Remark}
\newtheorem{corollary}[theorem]{Corollary}

\numberwithin{equation}{section}

\begin{document}

\baselineskip=15pt

\title[Rank one connections on abelian varieties]{Rank one
connections on abelian varieties}

\author[I. Biswas]{Indranil Biswas}

\address{School of Mathematics, Tata Institute of Fundamental
Research, Homi Bhabha Road, Bombay 400005, India}

\email{indranil@math.tifr.res.in}

\author[J. Hurtubise]{Jacques Hurtubise}

\address{Department of Mathematics, McGill University, Burnside
Hall, 805 Sherbrooke St. W., Montreal, Que. H3A 2K6, Canada}

\email{jacques.hurtubise@mcgill.ca}

\author[A. K. Raina]{A. K. Raina}

\address{Theoretical Physics Department, Tata Institute of
Fundamental Research, Homi Bhabha Road, Bombay 400005, India}

\email{raina@tifr.res.in}



\date{}

\begin{abstract}
Let $A$ be a complex abelian variety. The moduli space
${\mathcal M}_C$ of rank one algebraic
connections on $A$ is a principal bundle over the dual abelian
variety $A^\vee\,=\,\text{Pic}^0(A)$ for the group $H^0(A,\, \Omega^1_A)$. 
Take any
line bundle $L$ on $A^\vee$; let ${\mathcal C}(L)$ be the
algebraic principal $H^0(A^\vee, \Omega^1_{A^\vee})$--bundle 
over $A^\vee$ given by the sheaf of connections on $L$. The line 
bundle $L$ produces a homomorphism $H^0(A,\, \Omega^1_A)
\,\longrightarrow\, H^0(A^\vee,\, \Omega^1_{A^\vee})$. We prove that
${\mathcal C}(L)$ is isomorphic to the principal
$H^0(A^\vee, \Omega^1_{A^\vee})$--bundle obtained by extending the
structure group of the principal $H^0(A,\, \Omega^1_A)$--bundle
${\mathcal M}_C$ using this homomorphism given by $L$. We
compute the ring of algebraic functions on ${\mathcal C}(L)$. As
an application of the above result, we show that ${\mathcal M}_C$
does not admit any non-constant algebraic function, despite the
fact that it is biholomorphic to $({\mathbb C}^*)^{2\dim A}$
implying that it has many non-constant holomorphic functions.
\end{abstract}

\maketitle

\section{Introduction}

Let $A$ be a complex abelian variety of dimension $d_0$. Let $\mathcal M$
be the moduli space of rank one $\lambda$--connections on $A$ (the
definition of $\lambda$--connections is recalled in Section \ref{sec2}).
This $\mathcal M$ is a vector bundle over $A^\vee\,:=\,\text{Pic}^0(A)$, and
it fits in the following short exact sequence of vector
bundles on $A^\vee$
\begin{equation}\label{i1}
0\,\longrightarrow\, A^\vee\times H^0(A,\, \Omega^1_A)\,
\longrightarrow\,{\mathcal M}\,\longrightarrow\,
{\mathcal O}_{A^\vee}\,\longrightarrow\, 0\, .
\end{equation}

Take any line bundle $L$ on $A^\vee$. It produces a homomorphism
\begin{equation}\label{i3}
\phi^*_L\, :\, H^0(A,\, \Omega^1_A)\,\longrightarrow\,
H^0(A^\vee,\, \Omega^1_{A^\vee})\, ;
\end{equation}
this homomorphism is the pull back of differential forms by the map 
$A^\vee\, \longrightarrow\, A\,=\,\text{Pic}^0(A^\vee)$ defined by 
$x\, \longmapsto\,
L^*\otimes \tau^*_xL$, where $\tau_x$ is the translation automorphism
$z\, \longmapsto\, z+x$ of $A^\vee$. Let
\begin{equation}\label{i2}
0\,\longrightarrow\,A^\vee\times H^0(A^\vee,\, \Omega^1_{A^\vee})\,=\,
\Omega^1_{A^\vee}
\,\longrightarrow\,{\mathcal W}\,\longrightarrow
\,{\mathcal O}_{A^\vee}\,\longrightarrow\, 0
\end{equation}
be the push--forward of the exact sequence in \eqref{i1} using the
homomorphism $\phi^*_L$ in \eqref{i3}.

Let
$$
0\, \longrightarrow\, \Omega^1_{A^\vee}\, \longrightarrow\, \text{At}(L)^*
\, \longrightarrow\, {\mathcal O}_{A^\vee} \, \longrightarrow\, 0
$$
be the dual of the Atiyah exact sequence for $L$. We prove that it is related
to the one in \eqref{i2} in the following way (see Proposition \ref{prop1}):

\begin{proposition}\label{prop0}
There exists an isomorphism of vector bundles $\eta\, :\,{\mathcal W}
\, \longrightarrow\, {\rm At}(L)^*$ such that the following
diagram is commutative
$$
\begin{matrix}
0 &\longrightarrow & \Omega^1_{A^\vee}
&\longrightarrow & {\mathcal W}&\longrightarrow& {\mathcal O}_{A^\vee}
&\longrightarrow & 0\\
&& \Vert && ~\Big\downarrow \eta && ~\,~\Big\downarrow \text{-}1\cdot \\
0& \longrightarrow & \Omega^1_{A^\vee}& \longrightarrow & \text{At}(L)^*
& \longrightarrow & {\mathcal O}_{A^\vee} & \longrightarrow & 0
\end{matrix}
$$
\end{proposition}

Let ${\mathcal M}_C$ be the moduli space of rank one
algebraic connections on $A$.
It is a principal bundle over $A^\vee$ for the additive group
$H^0(A,\, \Omega^1_A)$. Let ${\mathcal C}(L)$ be the fiber bundle over
$A^\vee$ given by the sheaf of algebraic
connections on $L$; it is in fact a
principal $H^0(A,\, \Omega^1_{A^\vee})$--bundle over $A^\vee$.

We prove the following (see Corollary \ref{cor2}):

\begin{corollary}\label{cor0}
The principal $H^0(A,\, \Omega^1_{A^\vee})$--bundle
${\mathcal C}(L)$ over $A^\vee$ is isomorphic to the principal
$H^0(A,\, \Omega^1_{A^\vee})$--bundle obtained by extending the
structure group of the principal $H^0(A,\, \Omega^1_A)$--bundle
${\mathcal M}_C$ using the homomorphism $\phi^*_L$ in \eqref{i3}.
\end{corollary}

Let $V^*_L$ be the cokernel of the homomorphism $\phi^*_L$.
There is a projection
$$
\rho\, :\, {\mathcal C}(L)\,\longrightarrow\, V^*_L\, .
$$
We prove the following theorem (see Theorem \ref{thm1}):

\begin{theorem}\label{thm0}
All algebraic functions on ${\mathcal C}(L)$ factor through the
projection $\rho$. In particular, if $\phi^*_L$ is an isomorphism,
then there are no non-constant algebraic functions on ${\mathcal C}(L)$.
\end{theorem}

We note that if either $L$ is ample or $L^*$ is ample, then
$\phi^*_L$ is an isomorphism.

Using Theorem \ref{thm0} it can be shown that there are no non-constant
algebraic functions on the moduli space ${\mathcal M}_C$. (See
Proposition \ref{prop3}.)

We note that ${\mathcal M}_C$ is canonically biholomorphic to 
${\rm Hom}(\pi_1(A), {\mathbb C}^*)$; the biholomorphism sends
an integrable connection on $A$ to its monodromy. Therefore, 
${\mathcal M}_C$ is
biholomorphic to $({\mathbb C}^*)^{2d_0}$. Consequently,
${\mathcal M}_C$ admits non-constant holomorphic functions.

Proposition \ref{prop3} yields the following theorem on the
moduli space of rank one $\lambda$--connections (see Theorem 
\ref{thm2}):

\begin{theorem}\label{thm-1}
All algebraic functions on $\mathcal M$ factor through
the forgetful map
$$
{\mathcal M}\, \longrightarrow\,{\mathbb C}
$$
that sends any $\lambda$--connection $(L'\, ,c'\, ,D')$ to $c'$. 
\end{theorem}

\section{Moduli space of $\lambda$--connections}\label{sec2}

Let $A$ be an abelian variety, defined over $\mathbb C$,
of dimension $d_0$, with $d_0\, \geq\, 1$.
The identity element of $A$ will be denoted by $0$.
The dual abelian variety $\text{Pic}^0(A)$ will be denoted
by $A^\vee$.

Consider a triple $(L\, ,c\, ,D)$, where $L$ is
an algebraic line bundle over $A$, $c\,\in\, \mathbb C$,
and
$$
D\, :\, L\, \longrightarrow\, L\otimes \Omega^1_A
$$
is an algebraic differential operator, of order at most
one, satisfying the identity
$$
D(f_0s) \,=\, f_0\cdot D(s) + cs\otimes df_0\, ,
$$
where $f_0$ is a locally defined algebraic function, and $s$
is a locally defined algebraic section of $L$. It is easy to
see that the operator $D$ is integrable, meaning the composition
$$
D\circ D\, :\, L\, \longrightarrow\, L\otimes \Omega^2_A
$$
vanishes identically. If $c\,=\, 0$,
then $D$ is an algebraic one--form on $A$, so
$(L\, ,0\, ,D)$ is a Higgs line bundle. If
$c\,\not=\, 0$, then $D/c$ is an integrable algebraic connection
on $L$. Therefore, if $c\,\not=\, 0$, then $L\,\in\,
\text{Pic}^0(A)\,=\, A^\vee$. 

A $\lambda$--\textit{connection} on $A$ is a triple
$(L\, ,c\, ,D)$ of the above form such that
$L\,\in\,A^\vee$. As mentioned above, if $c\,\not=\, 0$,
then the condition on $L$ is automatically satisfied.

If $(L\, ,c_1\, ,D_1)$ and $(L\, ,c_1\, ,D_1)$ are
two $\lambda$--connections, then $(L\, ,c_1+c_2\, ,D_1+D_2)$
is also a $\lambda$--connection. Similarly, for any
$\alpha\, \in\, \mathbb C$, the triple
$(L\, ,\alpha\cdot c_1\, ,\alpha\cdot D_1)$ is a $\lambda$--connection.
Therefore, the space of all $\lambda$--connections on a fixed
line bundle is a complex vector space. The dimension of this
vector space is $d_0+1$.

Let $(L\, ,c\, ,D)$ be a $\lambda$--connection on
$A$. For any $\omega\,\in\,
H^0(A, \, \Omega^1_A)$, the triple $(L\, ,c\, ,D+\omega)$ is also a
$\lambda$--connection. Moreover, if $(L\, ,c\, ,D')$ is another
$\lambda$--connection, then there is a unique $\omega'
\,\in\, H^0(A, \, \Omega^1_A)$ such that $D'\,=\, D+\omega'$. In other
words, for fixed $L$ and $c$, the space of all $\lambda$--connections
of the form $(L\, ,c\, ,D'')$ is an affine space for $H^0(A,\,\Omega^1_A)$.

Let ${\mathcal M}$ be the moduli space of $\lambda$--connections
on $A$. The moduli space of $\lambda$--connections of arbitrary
but fixed rank was constructed in \cite{Si}. Here we
are considering only $\lambda$--connections of rank one. It is
substantially simpler
to construct the moduli space of $\lambda$--connections of rank one.

Let
\begin{equation}\label{e1}
f\, :\, {\mathcal M}\, \longrightarrow\, A^\vee
\end{equation}
be the projection defined by $(L\, ,c\, ,D)\,\longmapsto\, L$.
It was noted above that each fiber of $f$ is a vector space of
dimension $d_0+1$. The projection $f$ makes ${\mathcal M}$ the
total space of an algebraic vector bundle over $A^\vee$ of
rank $d_0+1$. This vector bundle will also be denoted by
$\mathcal M$. Let
\begin{equation}\label{e2}
q_0\, :\, {\mathcal M}\, \longrightarrow\, A^\vee\times \mathbb C
\end{equation}
be the projection defined by $(L\, ,c\, ,D)\,\longmapsto\, (L\, ,c)$.
Let
\begin{equation}\label{a}
p_{\mathbb C}\, :\, A^\vee\times {\mathbb C}\, \longrightarrow\,
\mathbb C
\end{equation}
be the natural projection.
The restriction of the morphism $p_{\mathbb C}\circ
q_0$ to a fiber of $f$ is linear projection (recall that
the fibers of $f$ are vector spaces). The
kernel of this linear projection is $H^0(A,\, \Omega^1_A)$.
Therefore, we
get a short exact sequence of vector bundles on $A^\vee$
\begin{equation}\label{e3}
0\,\longrightarrow\, A^\vee\times H^0(A,\, \Omega^1_A)\,
\stackrel{\iota}{\longrightarrow}\,{\mathcal M} 
\,\stackrel{q_0}{\longrightarrow}\, {\mathcal O}_{A^\vee}
\,=\,A^\vee\times{\mathbb C}\,\longrightarrow\, 0\, ,
\end{equation}
where $A^\vee\times H^0(A,\, \Omega^1_A)$ is the trivial vector bundle
over $A^\vee$ with fiber $H^0(A,\, \Omega^1_A)$. This trivial vector
bundle $A^\vee\times H^0(A,\, \Omega^1_A)$ will be denoted by ${\mathcal 
V}_0$.

Let
\begin{equation}\label{e4}
\delta \, \in\, H^1(A^\vee,\, {\mathcal V}_0)\,=\, H^1(A^\vee,\,
{\mathcal O}_{A^\vee})\otimes H^0(A,\, \Omega^1_A)\,=\,
H^1(A^\vee,\, {\mathcal O}_{A^\vee})\otimes H^0(A,\, TA)^*
\end{equation}
be the extension class for the short exact sequence in \eqref{e3}.

There is a canonical isomorphism
\begin{equation}\label{e4a}
H^0(A,\, TA)\, \longrightarrow\, H^1(A^\vee,\, {\mathcal O}_{A^\vee})
\end{equation}
(see \cite[p. 4, Theorem]{Mu}); note that since $A\,=\, (A^\vee)^\vee$,
and the homomorphism of evaluation of
sections at the origin $H^0(A,\, TA)\, \longrightarrow\,
T_0A$ is an isomorphism, the isomorphism in \eqref{e4a} follows from the
above mentioned result of \cite{Mu}. (We will recall a
construction of the isomorphism \eqref{e4a} in the proof
of Lemma \ref{lem1}.) The isomorphism in \eqref{e4a} 
gives an element
\begin{equation}\label{e5}
c_0\, \in\, H^1(A^\vee,\, {\mathcal O}_{A^\vee})\otimes H^0(A,\, TA)^*\, .
\end{equation}

\begin{lemma}\label{lem1}
The cohomology class $c_0$ in \eqref{e5} coincides with $-\delta$,
where $\delta$ is the class in \eqref{e4}.
\end{lemma}

\begin{proof}
Express $A$ as $V/\Lambda$, where $V$ is a complex vector space
of dimension $d_0$, and $\Lambda$ is a cocompact lattice in $V$.
Let $\overline{\Lambda}$ be the lattice in $\overline{V}^*$ given
by $\Lambda$. This lattice $\overline{\Lambda}$ is defined by the
condition that the natural pairing $\Lambda \times 
\overline{\Lambda}\,\longrightarrow\,
\mathbb C$ takes values in $\mathbb Z$.
The abelian variety $A^\vee$ is the quotient 
$\overline{V}^*/\overline{\Lambda}$.

We note that $H^0(A,\, TA)\,=\,
V$, because $H^0(A,\, TA) \,=\, T_0A \,=\, T_0V\,=\, V$
(the fiber $T_0A$ is identified
with $H^0(A,\, TA)$ by evaluating vector fields at $0$). So,
$H^0(A,\, \Omega^1_A)\,=\, V^*$.
The conjugate--linear homomorphism
$$
H^0(A,\, \Omega^1_A)\,\longrightarrow\, H^1(A,\, {\mathcal O}_A)
$$
defined by $\theta\, \longmapsto\, \overline{\theta}$ is an isomorphism
(the closed form $\overline{\theta}$ represents a class in the
Dolbeault cohomology). Hence
\begin{equation}\label{b1}
H^1(A,\, {\mathcal O}_A) \,=\, \overline{V}^*\, .
\end{equation}
Consequently,
\begin{equation}\label{i}
H^1(A^\vee,\, {\mathcal O}_{A^\vee}) \,=\, \overline{\overline{V}^*}^*
\,=\, V\, .
\end{equation}
The isomorphism in \eqref{e4a} is the identity map of $V$.

A line bundle $L$ over $A$ with $c_1(L)\,=\, 0$ admits a unique
unitary flat connection. This unitary flat connection on $L$ will be
denoted by $\nabla^L$. So we get a $C^\infty$ section of the vector
bundle $\mathcal M$
\begin{equation}\label{s}
s\, :\, A^\vee \, \longrightarrow\, \mathcal M
\end{equation}
defined by $L\, \longmapsto\, (L\, ,1\, ,\nabla^L)$.
Note that the map
$$
{\mathcal O}_{A^\vee} \,\longrightarrow\, {\mathcal M} 
$$
defined by $(L\, ,c)\, \longmapsto\, c\cdot s(L)$ gives a
$C^\infty$ splitting of the short exact sequence in \eqref{e3}.

Let $\overline{\partial}_{\mathcal M}\, :\, C^\infty(A^\vee,\,
{\mathcal M})
\, \longrightarrow\, C^\infty(A^\vee,\, {\mathcal M}\otimes
(T^{0,1}A^\vee)^*)$ be the Dolbeault operator for the
holomorphic vector bundle $\mathcal M$. Since $q_0$
in \eqref{e3} is holomorphic,
$$
(q_0\times\text{Id})(\overline{\partial}_{\mathcal M}(s))
\,=\, \overline{\partial}(q_0(s))\,=\, \overline{\partial}(1)
\,=\, 0\, .
$$
Consequently, $\overline{\partial}_{\mathcal M}(s)$
is a $(0,1)$--form on $A^\vee$ with values in the vector bundle
${\mathcal V}_0\,=\, A^\vee\times H^0(A,\, \Omega^1_A)$.
The extension class
$\delta$ in \eqref{e4} is represented by the ${\mathcal V}_0$--valued
$(0,1)$--form $\overline{\partial}_{\mathcal M}(s)$.

Consider the trivial $C^\infty$ line bundle $A\times \mathbb C$ over $A$.
The $(0,1)$--component of the de Rham differential, namely 
$\overline{\partial}$,
defines the trivial holomorphic structure on it. Any other Dolbeault
operator on this $C^\infty$ line bundle is of the form 
$\overline{\partial} +\omega'$, where $\omega'$ is a $(0,1)$--form 
with $\overline{\partial}\omega'\, =\, 0$. For
any $\omega\, \in\, H^0(A,\, \Omega^1_A)$, the $(0\, ,1)$--form
$\omega'\,:=\, \overline{\omega}$ satisfies this condition.
Therefore, we have a map
\begin{equation}\label{psi}
\psi\, :\, \overline{H^0(A,\, \Omega^1_A)}\, \longrightarrow\, A^\vee
\end{equation}
that sends any $\omega\, \in\, H^0(A,\, \Omega^1_A)$ to the holomorphic
line bundle defined by the Dolbeault operation $\overline{\partial} +
\overline{\omega}$ on the trivial line bundle. This $\psi$ is a
holomorphic covering map.

Consider the section $s$ constructed in \eqref{s} and
the map $\psi$ in \eqref{psi}.
For any $\omega\, \in\, H^0(A,\, \Omega^1_A)$, the holomorphic connection
$s(\psi(\omega))$ is $\overline{\partial} - \omega$ on the topologically
trivial line bundle $A\times \mathbb C$, where $\overline{\partial}$ 
is the $(1,0)$--component of the de Rham differential. (The 
corresponding unitary connection
is $d+\overline{\omega} -\omega$.) Since the extension class
$\delta$ in \eqref{e4} is represented by
$\overline{\partial}_{\mathcal M}(s)$, we conclude
that $\delta$ coincides with $-c_0$.
\end{proof}

Take any algebraic line bundle $L$ over $A^\vee$. Using $L$, we get
a homomorphism
\begin{equation}\label{pL}
\phi_L\, :\, H^0(A^\vee, \,TA^\vee)\, \longrightarrow\, H^0(A,\, TA)\, ;
\end{equation}
we will recall the construction of $\phi_L$. Consider
$$
c_1(L)\, \in\, H^{1,1}(A^\vee) \,=\, H^{1}(A^\vee,\,
\Omega^1_{A^\vee})\, .
$$
Using the obvious pairing of $TA^\vee$ with its 
dual $\Omega^1_{A^\vee}$, we get a homomorphism
$$
H^{1,1}(A^\vee)\otimes H^0(A^\vee,\, TA^\vee)\, \longrightarrow\,
H^1(A^\vee,\, {\mathcal O}_{A^\vee})\, .
$$
Using it, $c_1(L)$ gives a homomorphism $H^0(A^\vee,\, TA^\vee)\,
\longrightarrow\,H^1(A^\vee,\, {\mathcal O}_{A^\vee})$. The
homomorphism $\phi_L$ in \eqref{pL} coincides with this homomorphism
using the identification
$$
H^1(A^\vee,\, {\mathcal O}_{A^\vee})\,=\, H^0(A,\, TA)
$$
(see \eqref{i}).

Another construction of $\phi_L$ is as follows. For any $x\, \in\,
A^\vee$, let $\tau_x$ be the translation automorphism of $A^\vee$
defined by $z\, \longmapsto\, x+z$. The map
$$
\widetilde{\phi}_L\, :\, A^\vee\,\longrightarrow\, (A^\vee)^\vee \,=\, A
$$
defined by $x\, \longmapsto \, (\tau^*_x L)\otimes L^*$ is a homomorphism
of algebraic groups (see \cite[pp. 59--60, Corollary 4]{Mu}). The
homomorphism $\phi_L$ in \eqref{pL} is the corresponding 
homomorphism of Lie algebras.

Fix a line bundle $L$ on $A^\vee$. Let
\begin{equation}\label{aes}
0\, \longrightarrow\, {\mathcal O}_{A^\vee}\, \longrightarrow\, 
\text{At}(L)\, \longrightarrow\, TA^\vee \, \longrightarrow\, 0
\end{equation}
be the Atiyah exact sequence for $L$;
see \cite[p. 187, Theorem 1]{At} for the construction of Atiyah exact 
sequence. Let
\begin{equation}\label{e6}
0\, \longrightarrow\, \Omega^1_{A^\vee}\, \longrightarrow\, \text{At}(L)^*
\, \stackrel{p_1}{\longrightarrow}\, {\mathcal O}_{A^\vee} \, 
\longrightarrow\, 0
\end{equation}
be the dual of the exact sequence in \eqref{aes}.

Let
\begin{equation}\label{e7}
\phi^*_L\, :\, H^0(A,\, \Omega^1_A)\,\longrightarrow\,
H^0(A^\vee,\, \Omega^1_{A^\vee})
\end{equation}
be the dual of the homomorphism $\phi_L$ in \eqref{pL}.
Using $\phi^*_L$, we will construct a short exact sequence from
the one in \eqref{e3}.

Consider the direct sum of vector bundles
$(A^\vee\times H^0(A^\vee,\, \Omega^1_{A^\vee}))\oplus{\mathcal M}$
over $A^\vee$, where $A^\vee\times H^0(A^\vee,\, \Omega^1_{A^\vee})$
is the trivial vector bundle with fiber $H^0(A^\vee,\, 
\Omega^1_{A^\vee})$. Let
$$
A^\vee\times H^0(A,\, \Omega^1_A)\,\longrightarrow\, (A^\vee\times 
H^0(A^\vee,\, \Omega^1_{A^\vee}))\oplus {\mathcal M}
$$
be the embedding defined by $(x\, ,v)\, \longmapsto\,
((x\, ,-\phi^*_L(v))\, , \iota (x\, ,v))$, where $\iota$ is the
homomorphism in \eqref{e3} and $\phi^*_L$ is constructed
in \eqref{e7}. Let
\begin{equation}\label{cW}
{\mathcal W}\, :=\, ((A^\vee\times H^0(A^\vee,\, \Omega^1_{A^\vee}))\oplus 
{\mathcal M})/(A^\vee\times H^0(A,\, \Omega^1_A))
\end{equation}
be the quotient of this injective homomorphism. Note that
${\mathcal W}$ has a 
projection to ${\mathcal O}_{A^\vee}$ defined by
$(a\, ,b)\, \longmapsto\, q_0 (b)$, where $q_0$ is the projection
in \eqref{e3}, and $a\, \in\, A^\vee\times H^0(A^\vee,\, 
\Omega^1_{A^\vee})$. Using this projection, the vector bundle
${\mathcal W}$ fits in a short exact
sequence
\begin{equation}\label{e8}
0\,\longrightarrow\,A^\vee\times H^0(A^\vee,\, \Omega^1_{A^\vee})\,=\,
\Omega^1_{A^\vee}
\,\longrightarrow\,{\mathcal W}\,\stackrel{p_2}{\longrightarrow}
\,{\mathcal O}_{A^\vee}\,\longrightarrow\, 0\, .
\end{equation}

\begin{proposition}\label{prop1}
There exists an isomorphism of vector bundles $\eta\, :\,{\mathcal W}
\, \longrightarrow\, {\rm At}(L)^*$ such that the following
diagram is commutative
$$
\begin{matrix}
0 &\longrightarrow & \Omega^1_{A^\vee}
&\longrightarrow & {\mathcal W}&\longrightarrow& {\mathcal O}_{A^\vee}
&\longrightarrow & 0\\
&& \Vert && ~\Big\downarrow \eta && ~\,~\Big\downarrow \text{-}1\cdot \\
0& \longrightarrow & \Omega^1_{A^\vee}& \longrightarrow & \text{At}(L)^*
& \longrightarrow & {\mathcal O}_{A^\vee} & \longrightarrow & 0
\end{matrix}
$$
where and the top and bottom exact
sequences are those in \eqref{e8} and \eqref{e6} respectively.
\end{proposition}

\begin{proof}
The extension class in $H^1(A^\vee,\, \Omega^1_{A^\vee})$
for the Atiyah exact sequence for $L$ in \eqref{aes} coincides with
$c_1(L)$ \cite[p. 196, Proposition 12]{At}. Therefore, the extension 
class for the dual exact sequence in \eqref{e6} is also $c_1(L)$.

Using the homomorphism $\phi^*_L$ in \eqref{e7}, the cohomology class
$$
c_0\, \in\, H^1(A^\vee,\, {\mathcal O}_{A^\vee})\otimes H^0(A,\, 
\Omega^1_A)
$$
in \eqref{e5} produces a cohomology class
\begin{equation}\label{c0}
\widetilde{c}_0\,\in\, H^1(A^\vee,\, {\mathcal O}_{A^\vee})\otimes
H^0(A^\vee,\, \Omega^1_{A^\vee})\,=\,H^1(A^\vee,\, \Omega^1_{A^\vee})\, .
\end{equation}
Similarly, using $\phi^*_L$, the cohomology class
$\delta\, \in\, H^1(A^\vee,\, {\mathcal O}_{A^\vee})\otimes H^0(A,\,
\Omega^1_A)$ in \eqref{e4} produces a cohomology class
\begin{equation}\label{d0}
\widetilde{\delta}_0\,\in\, H^1(A^\vee,\, \Omega^1_{A^\vee})\, .
\end{equation}
{}From the constructions of the vector bundle $\mathcal W$ and the
sequence in \eqref{e8} it follows immediately that the extension
class for the short exact sequence in \eqref{e8} coincides with
$\widetilde{\delta}_0$.

{}From Lemma \ref{lem1} it follows that $-\widetilde{c}_0$ in \eqref{c0}
coincides with $\widetilde{\delta}_0$ in \eqref{d0}. Therefore,
the extension class for the short exact sequence in \eqref{e8} is
$-\widetilde{c}_0$.

The cohomology class $\widetilde{c}_0$ coincides with
$c_1(L)$; this follows from the expression of $c_1(L)$
in \cite[p. 18, Proposition]{Mu} and the computation
in \cite[p. 83]{Mu}. Since the extension classes
for the short exact sequences in \eqref{e8} and \eqref{e6}
are $-\widetilde{c}_0$ and $\widetilde{c}_0$ respectively, the
proof of the proposition is complete.
\end{proof}

\begin{remark}\label{rem1}
{\rm It should be emphasized that Proposition \ref{prop1} does
not imply that the isomorphism class of the short exact sequence
in \eqref{e6} is independent of $L$, because the homomorphism
$\phi^*_L$ depends on $L$. Note that the short exact sequence
in \eqref{e6} splits if and only if $L$ admits a holomorphic
connection \cite[p. 188, Definition]{At}. Hence the short
exact sequence in
\eqref{e6} splits if $L$ is topologically trivial, and it
does not split if $L$ is ample.}
\end{remark}

Consider the projections $p_1$ and $p_2$ in \eqref{e6} and
\eqref{e8} respectively. Note that
$$
p^{-1}_1(A^\vee\times\{1\}) ~\, ~\, \text{~and~} ~\, ~\,
p^{-1}_2(A^\vee\times\{1\})
$$
are principal bundles over $A^\vee$ for the
additive group $H^0(A^\vee,\, \Omega^1_{A^\vee})$.

We have the following corollary of Proposition
\ref{prop1}.

\begin{corollary}\label{cor1}
The two principal $H^0(A^\vee,\, \Omega^1_{A^\vee})$--bundles
$p^{-1}_1(A^\vee\times\{1\})$ and $p^{-1}_2(A^\vee\times\{1\})$
are isomorphic.
\end{corollary}

\begin{proof}
Consider the isomorphism $\eta$ in Proposition \ref{prop1}. The
restriction of $-\eta$ to $p^{-1}_2(A^\vee\times\{1\})$ is
the required isomorphism.
\end{proof}

Consider $q^{-1}_0(A^\vee\times\{1\})$, where $q_0$ is the projection
in \eqref{e3}. It is a principal bundle over $A^\vee$ for the
additive group $H^0(A,\, \Omega^1_A)$. From the constructions of the
vector bundle $\mathcal W$ and exact sequence in
\eqref{e8} it follows that the principal
$H^0(A^\vee,\, \Omega^1_{A^\vee})$--bundle
$p^{-1}_2(A^\vee\times\{1\})$ over $A^\vee$ is identified with the one 
obtained by
extending the structure group of the principal
$H^0(A,\, \Omega^1_A)$--bundle
$q^{-1}_0(A^\vee\times\{1\})$ using the homomorphism
$\phi^*_L$ in \eqref{e7}.

Therefore, the following is a reformulation of Corollary
\ref{cor1}.

\begin{corollary}\label{cor2}
The principal
$H^0(A^\vee,\, \Omega^1_{A^\vee})$--bundle
$p^{-1}_1(A^\vee\times\{1\})$ is isomorphic to the one obtained by
extending the structure group of the
principal $H^0(A,\, \Omega^1_A)$--bundle
$q^{-1}_0(A^\vee\times\{1\})$ using the homomorphism
$\phi^*_L$.
\end{corollary}

\section{Algebraic functions on sheaf of connections}\label{sec3}

Let $A$ be a complex abelian variety, with $\dim A \, >\, 0$.
Let $\xi$ be an algebraic line bundle over $A$. Let
\begin{equation}\label{e9}
0\, \longrightarrow\, \Omega^1_A\, \stackrel{\iota_1}{\longrightarrow}
\, \text{At}(\xi)^*
\, \stackrel{q}{\longrightarrow}\, {\mathcal O}_A \, 
\longrightarrow\, 0
\end{equation}
be the dual of the Atiyah exact sequence for $\xi$.
Let
\begin{equation}\label{f8}
q_1\, :\, \text{At}(\xi)^*\, \longrightarrow\, \mathbb C
\end{equation}
be the composition of $q$ and the projection from the
total space of ${\mathcal O}_A$ to $\mathbb C$. The sheaf of
sections of $\text{At}(\xi)^*$ is the sheaf of $\lambda$--connections
on the line bundle $\xi$. Define the fiber bundle over $A$
\begin{equation}\label{e10}
{\mathcal C}(\xi)\, :=\, q^{-1}_1(1)\, \subset\,
\text{At}(\xi)^*\, .
\end{equation}
The sheaf of holomorphic (respectively, algebraic) sections of
${\mathcal C}(\xi)$ is the sheaf of holomorphic (respectively, 
algebraic) connections on $\xi$. Using the homomorphism $\iota_1$
in \eqref{e9}, the fiber bundle ${\mathcal C}(\xi)$ is a principal
bundle over $A$ for the group $H^0(A, \,\Omega^1_A)$.

Let
\begin{equation}\label{e11}
{\mathbb P}\, :=\, P(\text{At}(\xi)^*) ~\,~\, \text{~and~}
~\,~\, {\mathbb P}_0\, :=\, P(\Omega^1_A)
\end{equation}
be the projective bundles over $A$ parametrizing lines in
$\text{At}(\xi)^*$ and $\Omega^1_A$ respectively. Using
$\iota_1$ in \eqref{e9}, the projective bundle ${\mathbb P}_0$
is a subbundle of ${\mathbb P}$, and
\begin{equation}\label{e12}
{\mathcal C}(\xi)\, =\, {\mathbb P}\setminus {\mathbb P}_0\, .
\end{equation}

Let ${\mathcal O}_{\mathbb P}(-1)\, \longrightarrow\,
{\mathbb P}$ be the tautological line bundle. Let
\begin{equation}\label{e13}
\sigma\, \in\, H^0({\mathbb P},\, {\mathcal O}_{\mathbb P}(1))
\end{equation}
be the section of ${\mathcal O}_{\mathbb P}(1)\,:=\,
{\mathcal O}_{\mathbb P}(-1)^*$ defined by
$\sigma(v)\,=\, q_1(v)\,\in\, {\mathbb C}$,
where $v\,\in\, {\mathcal O}_{\mathbb P}(-1)$, and $q_1$
is the projection in \eqref{f8}. The divisor ${\rm Div}(\sigma)$
coincides with ${\mathbb P}_0$ in \eqref{e11}.

Let ${\mathbb C}[{\mathcal C}(\xi)]$ be the space of
algebraic functions on the variety ${\mathcal C}(\xi)$. From
\eqref{e12} it follows that ${\mathbb C}[{\mathcal C}(\xi)]$
is identified with the direct limit
$$
\lim_{n\to \infty} H^0({\mathbb P},\, {\mathcal O}_{\mathbb P}(
n{\mathbb P}_0))\, .
$$
Since ${\rm Div}(\sigma)\,=\, {\mathbb P}_0$,
\begin{equation}\label{e14}
{\mathbb C}[{\mathcal C}(\xi)]\,=\,\lim_{n\to \infty} H^0(
{\mathbb P},\,{\mathcal O}_{\mathbb P}(n))\, ,
\end{equation}
where ${\mathcal O}_{\mathbb P}(n)\,=\,
{\mathcal O}_{\mathbb P}(1)^{\otimes n}$; the vector space
$H^0({\mathbb P},\,{\mathcal O}_{\mathbb P}(n))$ sits inside
$H^0({\mathbb P},\,{\mathcal O}_{\mathbb P}(n+1))$ using the
homomorphism $s\longmapsto\, s\otimes \sigma$.

Our aim in this section is to compute the algebra ${\mathbb C}
[{\mathcal C}(\xi)]$.

For any $x\,\in\, A$, let $\tau_x\, :\, A\, \longrightarrow
\,A$ be the automorphism defined by $z\,\longmapsto\,
z+x$. Let
$$
A_0\, \subset\, A
$$
be the closed subgroup defined by all $x$ such that
$\tau^*_x \xi\,=\, \xi$ (see \cite[pp. 59--60, Corollary 4]{Mu}).
Let
\begin{equation}\label{e15}
V_\xi\, :=\, \text{Lie}(A_0)\,=\, T_0A_0\, \subset
\, T_0A \,\subset\, H^0(A,\, TA)
\end{equation}
be the tangent space, where $0\,\in\, A$ is the
identity element. The subspace $V_\xi$ in \eqref{e15} defines
an algebraic foliation on $A$. This foliation will be denoted by
$\mathcal F$, so ${\mathcal F}\,=\, A\times V_\xi$.
The evaluation of sections of $\mathcal F$ at $0$ gives an
isomorphism
$$
H^0(A,\, {\mathcal F})\, \stackrel{\sim}{\longrightarrow}\,
V_\xi\, .
$$

Take any Zariski open subset $U\, \subset\, A$.
A \textit{partial connection} on $\xi_U\, :=\, \xi\vert_U$ in the
direction of $\mathcal F$ is an algebraic first order differential
operator
$$
D\, :\, \xi_U\, \longrightarrow\, \xi_U\otimes {\mathcal F}^*
$$
satisfying the identity
$$
D(f_0s) \,=\, f_0\cdot D(s) + \widetilde{df_0}\otimes s\, ,
$$
where $f_0$ is a locally defined function, $s$ is a locally
defined section of $\xi_U$, and $\widetilde{df_0}$ is the
projection of the $1$--form $df_0$ to ${\mathcal F}^*\vert_U$.

Let ${\mathcal C}_{\mathcal F}(\xi)$ be the sheaf of partial
connections on $\xi$ in the direction of $\mathcal F$.

We will give an alternative description of
${\mathcal C}_{\mathcal F}(\xi)$.

Consider the projection
$\iota^*_1\, :\, \text{At}(\xi)\, \longrightarrow\, TA$,
where $\iota_1$ is the homomorphism in \eqref{e9}. Define
$$
\text{At}_{\mathcal F}(\xi)\,=\, (\iota^*_1)^{-1}({\mathcal F})\, .
$$
We have a short exact sequence of vector bundles
\begin{equation}\label{d1}
0\, \longrightarrow\, {\mathcal F}^*\, \stackrel{\iota'_1}{\longrightarrow}
\, \text{At}_{\mathcal F}(\xi)^*
\, \stackrel{q'_1}{\longrightarrow}\, {\mathcal O}_A \,
\longrightarrow\, 0\, ;
\end{equation}
it coincides with \eqref{e9} if ${\mathcal F}\,=\, TA$.
The fiber bundle ${\mathcal C}_{\mathcal F}(\xi)\,\longrightarrow
\, A$ is identified
with the inverse image $(q'_1)^{-1}(A\times \{1\})$. We note that
${\mathcal C}_{\mathcal F}(\xi)\,=\, (q'_1)^{-1}(A\times \{1\})$
is a principal bundle over $A$ for the additive group
$H^0(A,\, \mathcal F)$.

We will show that $\xi$ admits a partial connection
in the direction of $\mathcal F$.

Let
\begin{equation}\label{e16}
\phi_\xi\, :\, H^0(A,\, TA)\, \longrightarrow\, H^0(A^\vee,
\, TA^\vee)
\end{equation}
be the homomorphism constructed as in \eqref{pL} using the
line bundle $\xi$.

\begin{lemma}\label{lem2}
The subspace $V_\xi$ in \eqref{e15} coincides with the kernel of
the homomorphism $\phi_\xi$ in \eqref{e16}.
\end{lemma}

\begin{proof}
This lemma follows immediately by comparing \cite[p. 18,
Proposition]{Mu} and \cite[p. 84, Proposition (iii)]{Mu}.
\end{proof}

\begin{lemma}\label{lem3}
The short exact sequence in \eqref{d1} splits.
\end{lemma}

\begin{proof}
We have a commutative diagram
\begin{equation}\label{et}
\begin{matrix}
0& \longrightarrow & \Omega^1_A & \longrightarrow
& \text{At}(\xi)^* & \longrightarrow & {\mathcal O}_A
& \longrightarrow & 0\\
&& \Big\downarrow && \Big\downarrow && \Vert\\
0& \longrightarrow & {\mathcal F}^* & \longrightarrow
& \text{At}_{\mathcal F}(\xi)^* & \longrightarrow & {\mathcal O}_A
& \longrightarrow & 0
\end{matrix}
\end{equation}
where the top and the bottom exact sequences are as in \eqref{e9} and
\eqref{d1} respectively. Therefore, the extension class for
the bottom exact sequence is given by the extension class for the
top exact sequence.

We have an isomorphism
\begin{equation}\label{re1}
H^{1}(A,\, \Omega^1_A) \, =\, \text{Hom}_{\mathbb C}(H^0(A,\, 
TA)\, , H^0(A^\vee,\, TA^\vee))\, ,
\end{equation}
because $H^1(A,\, ,{\mathcal O}_A)
\,=\, H^0(A^\vee,\, TA^\vee)$ (see \eqref{e4a}). By 
this isomorphism, the extension class in $H^{1}(A,\, \Omega^1_A)$ for
the short
exact sequence in \eqref{e9} is mapped to the homomorphism $\phi_\xi$ in 
\eqref{e16} (see the proof of Proposition \ref{prop1}).
We also have
$$
H^{1}(A,\, {\mathcal F}^*)\,=\, H^0(A,\, {\mathcal F}^*)\otimes
H^1(A,\, {\mathcal O}_A)\,=\, V^*_\xi\otimes H^0(A^\vee,\, TA^\vee)
$$
$$
=\,
\text{Hom}_{\mathbb C}(V_\xi\, , H^0(A^\vee,\, TA^\vee))\, ,
$$
where $V_\xi$ is defined in \eqref{e15}. Combining this
with the isomorphism in \eqref{re1} we get an isomorphism
$$
\text{Hom}_{\mathbb C}(H^1(A,\, \Omega^1_A)\, ,
H^1(A,\, {\mathcal F}^*))
$$
$$
\, =\,
\text{Hom}_{\mathbb C}(
\text{Hom}_{\mathbb C}(H^0(A,\, TA)\, ,H^0(A^\vee,\, TA^\vee))
\, ,\text{Hom}_{\mathbb C}(V_\xi\, , H^0(A^\vee,\, TA^\vee)))\, .
$$
This isomorphism takes the quotient map
$$
H^{1}(A,\, \Omega^1_A)\, \longrightarrow\, H^{1}(A,\, {\mathcal F}^*)
$$
to the natural quotient map
$$
\text{Hom}_{\mathbb C}(H^0(A,\, TA)\, ,H^0(A^\vee,\, A^\vee))
\, \longrightarrow\,\text{Hom}_{\mathbb C}(V_\xi\, ,
H^0(A^\vee,\, A^\vee))
$$
defined by the restriction of homomorphisms to the subspace $V_\xi$.

The extension class for the top exact sequence in \eqref{et}
is $c_1(\xi)$ \cite[p. 196, Proposition 12]{At}. The homomorphism 
$\phi_\xi$ in \eqref{e16} is sent to $c_1(\xi)$ by the
isomorphism in \eqref{re1} (see the proof of
Proposition \ref{prop1}). 

In view of the above description of
the quotient map
$$
H^{1}(A,\, \Omega^1_A)\, \longrightarrow\, H^{1}(A,\,
{\mathcal F}^*)\, ,
$$
{}from Lemma \ref{lem2} we now conclude that the extension
class in $H^{1}(A,\, {\mathcal F}^*)$ for the short exact 
sequence in \eqref{d1} vanishes. This completes the proof
of the lemma.
\end{proof}

Fix a splitting homomorphism
\begin{equation}\label{e17}
\gamma\, :\, \text{At}_{\mathcal F}(\xi)^*\, \longrightarrow\, V^*_\xi
\end{equation}
for the short exact sequence in \eqref{d1} (recall that ${\mathcal F}
\,=\, A\times V_\xi$); so $\gamma\circ \iota'_1 (x\, ,v)\,=\, v$,
where $x\,\in\, A$, $v\,\in\, V^*_\xi$, and $\iota'_1$ is the
homomorphism in \eqref{d1}.

Consider the dual homomorphism $\text{At}(\xi)^*\, \longrightarrow\,
\text{At}_{\mathcal F}(\xi)^*$ of the natural inclusion map
$\text{At}_{\mathcal F}(\xi)\,\hookrightarrow\, \text{At}(\xi)$. The
restriction of it to ${\mathcal C}(\xi)$ (defined in \eqref{e10})
coincides with the map that sends a locally defined connection on $\xi$
to the partial connection given by it. Let
\begin{equation}\label{e18}
\rho\, :\, {\mathcal C}(\xi)\, \longrightarrow\, V^*_\xi
\end{equation}
be the surjective map obtained by composing this restriction map
with $\gamma$ in \eqref{e17}.

For any nonnegative integer $d$, let ${\mathcal D}_d$ be the space
of algebraic functions on $V^*_\xi$ satisfying the condition that the 
order
of pole at infinity is at most $d$. So,
\begin{equation}\label{e19}
{\mathcal D}_d\,=\, \text{Sym}^d(V_\xi\oplus {\mathbb C})\, .
\end{equation}

In view of \eqref{e19}, the pull back of 
functions on $V^*_\xi$ by the map $\rho$ in \eqref{e18} yields an 
inclusion
\begin{equation}\label{e20}
\text{Sym}^d(V_\xi\oplus {\mathbb C})\, \subset\,
H^0({\mathbb P},\,{\mathcal O}_{\mathbb P}(d))
\end{equation}
(see \eqref{e14} for $H^0({\mathbb P},\,{\mathcal O}_{\mathbb P}(d))$).

\begin{proposition}\label{prop2}
For all $d\, \geq\, 1$,
$$\dim H^0({\mathbb P},\,{\mathcal O}_{\mathbb P}(d))\, \leq\,
\dim {\rm Sym}^d(V_\xi\oplus{\mathbb C})\, .$$
\end{proposition}

\begin{proof}
We first note that $H^0({\mathbb P},\,{\mathcal O}_{\mathbb P}(d))\,
=\, H^0(A,\, \text{Sym}^d({\rm At}(\xi)))$. Taking symmetric
power of the vector bundles in the Atiyah exact sequence
\begin{equation}\label{at}
0\,\longrightarrow\, {\mathcal O}_A \,\longrightarrow\, \text{At}(\xi)\,
\longrightarrow\, TA \,\longrightarrow\, 0\, ,
\end{equation}
we have a short exact sequence of vector bundles
\begin{equation}\label{e21}
0\,\longrightarrow\, \text{Sym}^n({\rm At}(\xi)) \,\longrightarrow\,
\text{Sym}^{n+1}({\rm At}(\xi))\, \longrightarrow\, \text{Sym}^{n+1}(TA)
\,\longrightarrow\, 0
\end{equation}
for each $n\, \geq \, 0$ (by convention, $\text{Sym}^0$ of a
nonzero vector space is $\mathbb C$). Let
\begin{equation}\label{e22}
0\,\longrightarrow\, H^0(A,\, \text{Sym}^n({\rm At}(\xi))) 
\,\longrightarrow\, H^0(A,\, \text{Sym}^{n+1}({\rm At}(\xi)))
\end{equation}
$$
\longrightarrow\, H^0(A,\, \text{Sym}^{n+1}(TA))\,
\stackrel{\beta'}{\longrightarrow}\, H^1(A,\, \text{Sym}^n({\rm At}(\xi)))
\,\longrightarrow\,
$$
be the long exact sequence of cohomologies for \eqref{e21}.

The projection $\text{Sym}^{n}({\rm At}(\xi))\, \longrightarrow\,
\text{Sym}^{n}(TA)$ in \eqref{e21} gives a homomorphism
$$
H^1(A,\, \text{Sym}^n({\rm At}(\xi)))
\,\longrightarrow\, H^1(A,\, \text{Sym}^{n}(TA))\, .
$$
Let
\begin{equation}\label{e23}
\beta_{n+1}\, :\, H^0(A,\, \text{Sym}^{n+1}(TA))\,\longrightarrow\,
H^1(A,\, \text{Sym}^{n}(TA))
\end{equation}
be the composition of this homomorphism and $\beta'$ in \eqref{e22}.

{}From \eqref{e22} it follows that for all $d\, \geq\, 1$,
\begin{equation}\label{e24}
\dim H^0(A,\, \text{Sym}^d({\rm At}(\xi))) \,\leq\,
1+ \sum_{i=1}^d\dim {\rm kernel}(\beta_i)\, ,
\end{equation}
where $\beta_i$ is constructed in \eqref{e23}.

We will denote the vector space $H^0(A,\, TA)$ by $W$. Then,
$$
H^0(A,\, \text{Sym}^{n+1}(TA))\, =\, \text{Sym}^{d+1}(W)\, .
$$
Also,
$$
H^1(A,\, \text{Sym}^{n}(TA))\,=\, H^1(A,\, {\mathcal O}_A)\otimes
H^0(A,\, \text{Sym}^{n}(TA))\, =\, \overline{W}^*\otimes
\text{Sym}^n(W)
$$
(see \eqref{b1}), and
$$
H^1(A,\, \Omega^1_A)\,=\, H^1(A,\, {\mathcal O}_A)\otimes
H^0(A,\, \Omega^1_A)\,=\, \overline{W}^*\otimes W^*\,=\,
\text{Hom}_{\mathbb C}(W\, ,\overline{W}^*)\, .
$$

The homomorphism
$$
\beta_{i+1}\, :\, H^0(A,\, \text{Sym}^{i+1}(TA))\,=\,\text{Sym}^{i+1}
(W)\,\longrightarrow\, H^1(A,\, \text{Sym}^{i}(TA))\,=\,
\overline{W}^*\otimes \text{Sym}^i(W)
$$
in \eqref{e23} is given by
$$
c_1(\xi)\,\in\,H^1(A,\, \Omega^1_A)\,=\,
\text{Hom}_{\mathbb C}(W\, ,\overline{W}^*)
$$
using the obvious pairing of $W$ with $W^*$. Indeed, this
follows from the fact that the extension class for the
short exact sequence in \eqref{at} is given by $c_1(\xi)$
using the pairing between $W$ and $W^*$. In view of Lemma
\ref{lem2}, we now conclude that
\begin{equation}\label{e25}
\text{kernel}(\beta_{i+1})\,=\, \text{Sym}^{i+1}(V_\xi)\,\subset
\, \text{Sym}^{i+1}(V)\, .
\end{equation}
The proposition follows from \eqref{e24} and \eqref{e25}.
\end{proof}

\begin{theorem}\label{thm1}
For all $d\, \geq\, 1$,
$${\rm Sym}^d(V_\xi\oplus {\mathbb C})\, =\,
H^0({\mathbb P},\,{\mathcal O}_{\mathbb P}(d))\, .$$
All algebraic functions on ${\mathcal C}(\xi)$ factor
through the surjective map $\rho$ in \eqref{e18}.
\end{theorem}

\begin{proof}
In view of Proposition \ref{prop2}, the inclusion homomorphism
$$
\text{Sym}^d(V_\xi\oplus {\mathbb C})\, \hookrightarrow\,
H^0({\mathbb P},\,{\mathcal O}_{\mathbb P}(d))
$$
in \eqref{e20} is an isomorphism.

Using \eqref{e14} together with the above isomorphism, it
follows that all algebraic functions on ${\mathcal C}(\xi)$
are of the form $g\circ\rho$, where $g$ is an algebraic
function on $V^*_\xi$.
\end{proof}

\section{Algebraic functions on the moduli space
of $\lambda$--connections}

In this section we will use the set--up and notation of Section 
\ref{sec2}.

Let ${\mathcal M}_C$ be the moduli space of algebraic connections
of rank one on $A$. So,
\begin{equation}\label{pq}
{\mathcal M}_C\,=\, q^{-1}_0(A\times \{1\})\, \subset\,
{\mathcal M}\, ,
\end{equation}
where $q_0$ is defined in \eqref{e2}.

\begin{proposition}\label{prop3}
There are no non-constant algebraic functions on ${\mathcal M}_C$.
\end{proposition}

\begin{proof}
Fix a line bundle $L\, \longrightarrow\, A^\vee$ such that the homomorphism
$\phi_L$ in \eqref{pL} is an isomorphism. For example, if $L$ is ample,
then $\phi_L$ is an isomorphism (this follows from the fact that
statement (3) in ``Theorem of Lefschetz'' in \cite[pp. 29--30]{Mu}
implies statement (2) in the same theorem). Since $\phi_L$ is
an isomorphism,
the homomorphism $\phi^*_L$ in \eqref{e7} is also an isomorphism.

Consider the vector bundle $\mathcal W$ constructed in \eqref{cW}. Since 
$\phi^*_L$ is an 
isomorphism, the map of fiber bundles
$$
\psi\, :\, {\mathcal M}\, \longrightarrow\, {\mathcal W}
$$
that sends any $z\, \in\, \mathcal M$ to the equivalence class of
$(0\, ,z)\, \in\, (A^\vee\times H^0(A^\vee,\, \Omega^1_{A^\vee}))\oplus 
{\mathcal M}$ is an isomorphism. Therefore,
$$
-\eta\circ\psi\, :\, {\mathcal M}\, \longrightarrow\,\text{At}(L^*)
$$
is an isomorphism, where $\eta$ is the isomorphism
in Proposition \ref{prop1}. Furthermore,
$$
p_1\circ(-\eta\circ\psi)\,=\, q_0\, ,
$$
where $q_0$ and $p_1$ are the homomorphisms in \eqref{e2} and
\eqref{e6} respectively. Now from the commutativity of the
diagram of maps in Proposition \ref{prop1} we have
\begin{equation}\label{g1}
p^{-1}_1(A^\vee\times\{1\})\,=\, q^{-1}_0(A\times \{1\})\,=\,
{\mathcal M}_C\, .
\end{equation}

Since the homomorphism $\phi_L$ is an isomorphism, from Theorem
\ref{thm1} we know that all algebraic functions on 
$p^{-1}_1(A^\vee\times\{1\})$ are constant functions. Hence from
\eqref{g1} it follows that there are no non-constant
functions on ${\mathcal M}_C$.
\end{proof}

As mentioned in the introduction, the moduli space
${\mathcal M}_C$ is biholomorphic to $${\rm Hom}(\pi_1(A), {\mathbb 
C}^*)\, \cong\, ({\mathbb C}^*)^{2d_0}$$ by sending an integrable
connection to its monodromy representation. Since the variety
$({\mathbb C}^*)^{2d_0}$ has nonconstant algebraic functions,
we conclude that ${\mathcal M}_C$ is not algebraically
isomorphic to $({\mathbb C}^*)^{2d_0}$. (See
\cite[p. 101, Remark]{Se} for an example of
similar phenomenon.)

\begin{theorem}\label{thm2}
All algebraic functions on $\mathcal M$ factor through the surjective
map
$$
p_{\mathbb C}\circ q_0\, :\, {\mathcal M}\, \longrightarrow\,{\mathbb C}\, ,
$$
where $q_0$ and $p_{\mathbb C}$ are defined in \eqref{e2} and \eqref{a}
respectively.
\end{theorem}

\begin{proof}
For any nonzero number $t\, \in\, \mathbb C$, the inverse image
$(p_{\mathbb C}\circ q_0)^{-1}(t)$ is isomorphic to ${\mathcal M}_C$.
Hence the theorem follows from Proposition \ref{prop3}.
\end{proof}

\medskip
\noindent
\textbf{Acknowledgements.}\, The first named author thanks
McGill University for hospitality while the work was carried out.


\end{document}